\documentclass[leqno,12pt]{article}

\usepackage{latexsym}
 \usepackage{graphicx}
\usepackage{amsmath}
\usepackage{amssymb}
\usepackage{algorithm}
\usepackage{algorithmic}
\usepackage{cite}
\usepackage{color}

\usepackage[margin=1.3in, top=1.3in, bottom=1.3in]{geometry}

\newcommand{\nc}{\newcommand}
\nc{\nt}{\newtheorem}
\nt{thm}{Theorem}[section]
\nt{cor}[thm]{Corollary}
\nt{prop}[thm]{Proposition}
\nt{lem}[thm]{Lemma}
\nt{defn}[thm]{Definition}
\nt{rem}[thm]{Remark}
\nt{exa}[thm]{Example}
\nt{ass}[thm]{Assumption}
\nt{alg}[thm]{Algorithm}
\nc{\ip}[2]{\mbox{$\langle #1,#2 \rangle$}}
\nc{\pf}{\noindent{\bf Proof\ \ }}
\nc{\finpf}{\hfill{$\Box$}\linespace}
\nc{\linespace}{\vspace{\baselineskip} \noindent}
\nc{\R}{{\bf R}}
\nc{\bR}{{\overline{\R}}}
\nc{\U}{{\bf U}}
\nc{\V}{{\bf V}}
\nc{\W}{{\bf W}}
\nc{\X}{{\bf X}}
\nc{\Y}{{\bf Y}}
\nc{\Z}{{\bf Z}}
\nc{\Rn}{{\bf R}^n}
\nc{\bx}{\bar{x}}
\nc{\by}{\bar{y}}
\nc{\inT}{\mbox{\rm int}\,}
\nc{\cl}{\mbox{\rm cl}\,}
\nc{\gph}{\mbox{\rm gph}\,}

\def\tto{\;{\lower 1pt \hbox{$\rightarrow$}}\kern -12pt
           \hbox{\raise 2.8pt \hbox{$\rightarrow$}}\;}
\newenvironment{myequation}{\setcounter{equation}{\value{thm}}
   \begin{equation}}{\addtocounter{thm}{1}\end{equation}}

\nc{\bmye}{\begin{myequation}}
\nc{\emye}{\end{myequation}}

\begin{document}
\title{
Partial smoothness and constant rank
}
\author{A.S. Lewis\thanks{ORIE, Cornell University, Ithaca, NY 14853, U.S.A.
\texttt{people.orie.cornell.edu/aslewis}.
Research supported in part by National Science Foundation Grant DMS-1208338 and by the US-Israel Binational Science Foundation Grant 2008261.}
\and
Jingwei Liang\thanks{DAMTP, University of Cambridge, U.K. 
\texttt{jl993@cam.ac.uk}}
}
\date{\today}
\maketitle

\begin{abstract}
The idea of partial smoothness in optimization blends certain smooth and nonsmooth properties of feasible regions and objective functions.  As a consequence, the standard first-order conditions guarantee that diverse iterative algorithms (and post-optimality analyses) identify active structure or constraints.  However, by instead focusing directly on the first-order conditions, the formal concept of partial smoothness simplifies dramatically:  in basic differential geometric language, it is just a constant-rank condition.  In this view, partial smoothness extends to more general mappings, such as saddlepoint operators underlying primal-dual splitting algorithms.
\end{abstract}
\medskip

\noindent{\bf Key words:} partial smoothness, active set identification, nonsmooth optimization, subdifferential, primal-dual splitting
\medskip

\noindent{\bf AMS 2000 Subject Classification:} \ldots
\medskip

\section{Introduction}
A variety of optimization algorithms, ranging from classical active set methods to contemporary first-order algorithms for machine learning and high-dimensional statistics, exhibit ``identification'' properties.  Iterates in the underlying Euclidean space $\U$ converging to an optimal solution $\bar u$ eventually lie in a subset $M \subset \U$ capturing the structure of the optimal solution.  In traditional nonlinear programming, $M$ might be the ``identifiable surface'' \cite{Wright} of the feasible region defined by the constraints active at optimality;  in machine learning applications, $M$ might consist of vectors with a certain sparsity pattern \cite{prx_lin}.

A simple but quite extensive model of this phenomenon, following the philosophy of \cite{Wright}, is as follows.  We consider minimizing a lower semicontinuous objective function $f \colon \U \to \overline\R$ (convex, for now), and assume that the set $M$ of interest is a smooth surface, or more precisely a manifold around $\bar u$, meaning that locally it consists of solutions of a system of \mbox{$C^{(2)}$-smooth} equations with linearly independent gradients.  Identification amounts to the property
\bmye \label{identify}
v_k \in \partial f(u_k), \quad u_k \to \bar u, \quad v_k \to 0
\quad \Rightarrow \quad u_k \in M ~~ \mbox{eventually},
\emye
where $\partial f$ denotes the classical subdifferential operator.  Earlier versions of this identifiability idea include \cite{Calamai-More87,Dunn87,Burke-More88,Burke90,Al-Khayyal-Kyparisis91,Ferris91,Flam92}.

Closely related to the identification property (\ref{identify}) is the idea that the function $f$ is {\em partly smooth} at the point $\bar u$ relative to the manifold $M$.  This property combines smoothness conditions on $f$ when restricted to $M$ with a sharpness property of 
$f$ in directions normal to $M$.  More precisely, around the point $\bar u$ the restrictions of the function $f$ and its subdifferential $\partial f$ to the manifold $M$ should be $C^{(2)}$-smooth and continuous respectively, and the affine span of $\partial f(\bar u)$ should be a translate of the normal space to $M$ at $\bar u$.  This property, along with the nondegeneracy assumption that zero lies in the relative interior of $\partial f(\bar u)$, together suffice to ensure identifiability (\ref{identify}), as shown in \cite[Thm 4.10]{shanshan}.  

As a simple example, the function $f$ on the space $\R^2$ defined by $f(x,y) = |x| + y^2$ is partly smooth at its minimizer $(0,0)$ relative to the manifold $\{0\} \times \R$, and zero lies in the relative interior of the subdifferential $\partial f(0,0) = [-1,1] \times \{0\}$.  Hence the identifiability property (\ref{identify}) holds, as is easy to verify directly.

The terminology and original definition of partly smooth sets and functions originated in \cite{Lewis-active}.  A closely related thread of research, known as ``${\mathcal VU}$ theory'', originated with \cite{LS97}, and includes \cite{LOS,MC02,MC05,MC04,HS09,MC03}.  Inevitably, it seems, the formal definition of partly smooth sets and functions, and their ${\mathcal VU}$ analogues, are rather involved.  The definition of an identifiable surface in \cite{Wright} is not simple either, despite the transparency of the identifiability property (\ref{identify}).

As an approach to identifiability, considering partly smooth functions seems roundabout:  our aim, the property (\ref{identify}), involves only the subdifferential operator 
$\partial f$, and not the underlying function $f$.  It turns out that we can indeed characterize identifiability more naturally through a simple and fundamental property of the underlying operator $\partial f$.  Simply put, if the graph of the operator 
(in the product space $\U \times \U$) is a smooth manifold around the point $(\bar u,0)$, 
and the canonical projection of nearby points $(u,v)$ in the graph to $u \in \U$ is {\em constant rank}  (meaning that the projected tangent spaces at those points have constant dimension), then the identifiability property (\ref{identify}) follows.

In summary, the notion of partial smoothness, and the closely related idea of identifiability, are in essence constant-rank properties.  This perspective not only clarifies our understanding of these powerful tools, but broadens their potential application beyond the basic optimality condition $0 \in \partial f(\bar u)$ to more general variational conditions.  As an example, we consider the saddlepoint optimality conditions associated with primal-dual splitting methods like the Chambolle-Pock algorithm \cite{chambolle-pock}.

\section{Manifolds} \label{Manifolds}
We begin a more formal development by summarizing some elementary ideas about manifolds.  Given a Euclidean space $\U$, we consider a set $M \subset \U$ that has the structure of a smooth manifold locally, around a point $\bar u \in M$.  By ``smooth'', we mean $C^{(1)}$-smooth, unless we state otherwise.  We can consider such sets $M$ using ``local coordinates'', as follows.  

We denote the open ball of radius $\delta > 0$ around the point $\bar u$ by $B_{\delta}(\bar u)$.
In elementary language, 
$M$ is a {\em smooth manifold around} $\bar u$ when there exists a Euclidean space $\W$ and a map 
$H \colon \W \to \U$ that is smooth around $0$, with the derivative
$\nabla H(0) \colon \W \to \U$ injective and $H(0) = \bar u$, and such that, for all small $\delta > 0$,
\[
M = H\big(B_{\delta}(\bar u)\big) ~~\mbox{around} ~\bar u.
\]
More formally \cite[Chapter 8]{Lee}, some open neighborhood of $\bar u$ in $M$ is an {\em embedded submanifold} of $\U$.
Any small vector $w \in \W$ constitutes the {\em local coordinates centered around} $\bar u$ for the point 
$H(w) \in M$.  The {\em tangent space} at such a point is given simply by 
\[
T_M\big(H(w)\big) = \mbox{Range}\big(\nabla H(w)\big).
\]
Its dimension (the {\em dimension} of $M$ around $\bar u$) is a constant, namely $\dim \W$.  The {\em normal space} is the orthogonal complement:
\[
N_M\big(H(w)\big) = \mbox{Null}\big(\nabla H(w)^*\big).
\]

Given another Euclidean space $\V$, a map $F \colon M \to \V$ is {\em smooth around} $\bar u$ when there exists a map $G \colon \U \to \V$ that is smooth around $\bar u$ and agrees with $F$ on a neighborhood of $\bar u$ in $M$.  In that case, the {\em rank} of $F$ at $\bar u$ is 
$\dim\!\big(\nabla G(\bar u)T_M(\bar u)\big)$.  Equivalently, $F$ is smooth around $\bar u$ when the composition $F \circ H$ is smooth around $0$, and its rank at $\bar u$ is then rank of the derivative $\nabla(F \circ H)(0) \colon \W \to \V$ as a linear map.

The map $H$ defines a diffeomorphism from the open ball $B_{\delta}(0) \subset \W$ (for small $\delta > 0$) to an open neighborhood of the point $\bar u$ in the manifold $M$.  We can describe the inverse of this diffeomorphism via a map $G \colon \U \to \W$, smooth around the point $\bar u$, and satisfying 
\bmye \label{inverse}
G\big( H(w)\big) = w~~ \mbox{for all small vectors $w \in \W$}.
\emye
The restriction $G|_M$, around $\bar u$, is the inverse of the diffeomorphism $H$.

Adopting a dual approach, we can equivalently define a set $M \subset \U$ to be a smooth manifold around a point $\bar u$ when there exists a Euclidean space $\X$ and a map $P \colon \U \to \X$ that is smooth around $\bar u$, with the derivative $\nabla P(\bar u) \colon \U \to \X$ surjective and $P(\bar u) = 0$, and such that
\[
M ~=~ P^{-1}(0) ~=~ \{ u \in \U : P(u) =0 \} ~~\mbox{around} ~\bar u.
\]
Then the tangent and normal spaces are given by
\begin{eqnarray*}
T_M(u) &=& \mbox{Null}(\nabla P(u))  \\
N_M(u) &=& \mbox{Range}(\nabla P(u)^*)
\end{eqnarray*}
at all points $u \in M$ near $\bar u$.  The normal space has the same dimension as $\X$.

We can naturally decompose the space $U$ as a direct sum:
\[
\U = T_M(\bar u) \oplus N_M(\bar u).
\]
With this decomposition, the two derivatives $\nabla H(0) \colon \W \to \U$ and 
$\nabla P(\bar u) \colon \U \to \X$ are given by
\begin{eqnarray*}
\nabla H(0)w &=& (Dw,0) \\
\nabla P(\bar u)(r,s) &=& Es 
\end{eqnarray*}
for some invertible linear maps $D \colon \W \to T_M(\bar u)$ and $E \colon N_M(\bar u) \to \X$.  Furthermore, the derivative $\nabla G(\bar u) \colon \U \to \W$, restricted to $T_M(\bar u)$, is just the inverse map $D^{-1}$.

\section{Partly smooth mappings}
We consider the canonical projection
$\mbox{proj} \colon \U \times \V \to \U$ defined by $\mbox{proj}(u,v) = u$.   

\begin{defn}[Partly smooth mappings]
{\rm
A set-valued mapping $\Phi \colon \U \tto \V$ is called {\em partly smooth at} a point $\bar u \in \U$ {\em for} a value $\bar v \in \Phi(\bar u)$ when the graph $\gph \Phi$ is a smooth manifold around $(\bar u, \bar v)$ and the projection $\mbox{proj}$ restricted to $\gph \Phi$ has constant rank around 
$(\bar u, \bar v)$.  The {\em dimension} of $\Phi$ at $\bar u$ for $\bar v$ is then just the dimension of its graph around $(\bar u,\bar v)$.
}
\end{defn}

\noindent
{\bf Note.}  An example is when the inverse mapping $\Phi^{-1} \colon \V \tto \U$ is locally single-valued, smooth and constant-rank around $\bar v$ for $\bar u$.  In this case, $\Phi$ is in particular ``strongly regular'' at $\bar u$ for $\bar v$.
\bigskip

By definition, the constant rank condition means that the subspace
\[
\mbox{proj} \big(T_{\mbox{\scriptsize \gph}\,\Phi}(u,v) \big)
\]
and its orthogonal complement (called, in variational analysis, the {\em coderivative} of the mapping 
$\Phi$)
\[
D^*\Phi(u,v)(0) ~=~ \big\{ w \in \U : (w,0) \in N_{\mbox{\scriptsize \gph}\,\Phi}(u,v) \big\},
\]
or equivalently, the subspace
\[
N_{\mbox{\scriptsize \gph}\,\Phi}(u,v) \cap (\U \times \{0\})
\]
all have constant dimension for points $(u,v)$ near $(\bar u,\bar v)$.  

Consider, for example, the set-valued mapping $\Phi \colon \R \tto \R$ defined by
\[
\Phi(u) = 
\left\{
\begin{array}{ll}
\{ \pm \sqrt{u} \} & (u \ge 0) \\
\emptyset & (u < 0).
\end{array}
\right.
\]
The graph is of $\Phi$ is the manifold $\{(u,v) \in \R^2 : u = v^2 \}$.  However, $\Phi$ is not partly smooth at $0$ for $0$, because the projection $\mbox{proj}$ restricted to $\gph\Phi$ has rank zero at the point $(0,0)$ but rank one at all nearby points.

\begin{prop} \label{active}
If a set-valued mapping $\Phi \colon \U \tto \V$ is partly smooth at a point $\bar u \in \U$ for a value 
$\bar v \in \Phi(\bar u)$, then there exists a set 
$M \subset \U$, uniquely defined around $\bar u$, that is a smooth manifold around $\bar u$, and satisfies
\[
M ~=~ \big\{ u \in B_\epsilon(\bar u) : \exists v \in \Phi(u) \cap B_\epsilon(\bar v) \big\}~~ \mbox{around}~ \bar u,
\]
for all small $\epsilon > 0$.  We call any such set $M$ the {\bf active manifold}.
\end{prop}

\pf
For any small $\epsilon > 0$, the set
\[
G_{\epsilon} ~=~ \gph \Phi \cap \big( B_\epsilon(\bar u) \times B_{\epsilon}(\bar v) \big)
\]
is a manifold, and the projection $\mbox{proj}$ restricted to $G_{\epsilon}$ is a constant-rank map.  By the Constant Rank Theorem, the resulting image
\[
M_\epsilon ~=~ 
\big\{ u \in B_\epsilon(\bar u) : \exists v \in \Phi(u) \cap B_\epsilon(\bar v) \big\}
\]
is a manifold of dimension $\dim \mbox{proj}\,T_{G_{\epsilon}}(\bar u,\bar v)$.  This dimension is constant, for small $\epsilon > 0$, since the tangent space satisfies 
$T_{G_{\epsilon}}(\bar u,\bar v) = T_{\mbox{\scriptsize gph\,} \Phi}(\bar u,\bar v)$.  For any 
$\epsilon' \in (0,\epsilon)$, we know $M_{\epsilon'} \subset M_{\epsilon}$, but these sets are manifolds around $\bar u$ of the same dimension, so must be identical around $\bar u$.
\finpf

We use the following definition \cite{ident}.

\begin{defn}
{\rm
A set $M \subset \U$ is {\em identifiable} for a set-valued mapping $\Phi \colon \U \tto \V$ at a point 
$\bar u \in \U$ for a value $\bar v \in \Phi(\bar u)$ when $\gph \Phi \subset M \times \V$ around the point 
$(\bar u, \bar v)$.
}
\end{defn}

The following proposition is then immediate.

\begin{prop} \label{inclusion}
If a set-valued mapping $\Phi \colon \U \tto \V$ is partly smooth at a point $\bar u \in \U$ for a value 
$\bar v \in \Phi(\bar u)$, then the active manifold is an identifiable set.
\end{prop}
In fact, as we see shortly, the active manifold is a locally minimal identifiable set.

\section{Representations of partly smooth mappings}

The following result gives a representation of a partly smooth mapping using local coordinates.

\begin{thm}[Coordinate representation] \label{coordinate}
A set-valued mapping $\Phi \colon \U \tto \V$ is partly smooth at a point 
$\bar u \in \U$ for a value $\bar v \in \Phi(\bar u)$ if and only if it has a local representation of the following form:  there exist Euclidean spaces $\W$ and $\Z$, maps $H \colon \W \to \U$, smooth around $0$ with $H(0) = \bar u$ and $\nabla H(0)$ injective, and $G \colon \W \times \Z \to \V$, smooth around $(0,0)$ with $G(0,0) = \bar v$, such that, 
\bmye \label{regularity}
w \in \W,~ z \in \Z,~ \nabla H(0)w = 0~ \mbox{and}~ \nabla G(0,0)(w,z) = 0
~~\Rightarrow~~ w=0 ~\mbox{and}~ z=0,
\emye
and for all small $\delta > 0$,
\bmye \label{graph}
\gph \Phi ~=~ 
\big\{ \big(H(w),G(w,z)\big) :  w \in B_{\delta}(0),~  z \in B_{\delta}(0) \big\}
~~\mbox{around}~ (\bar u,\bar v).
\emye
In this case, the dimension of $\Phi$ at $\bar u$ for $\bar v$ is $\dim\W + \dim\Z$, and 
the active manifold is $H\big( B_{\delta}(0) \big)$ around $\bar u$,
providing $\delta > 0$ is sufficiently small.
\end{thm}

\pf
Assuming the local representation, we first prove that $\Phi$ is partly smooth at $\bar u$ for $\bar v$.  Consider the map $P \colon  \W \times \Z \to \U \times \V$ defined by $P(w,z) = \big(H(w),G(w,z)\big)$ for $w \in \W$ and $z \in \Z$.  This map is smooth around the point $(0,0)$, with derivative
\[
\nabla P(w,z)(r,s) = \big(\nabla H(w)r , \nabla G(w,z)(r,s)\big),
\]
for all small $w \in \W$ and $z \in \Z$, and vectors $r \in \W$ and $s \in \Z$.  By assumption, the derivative $\nabla P(0,0)$ is injective, so $\gph \Phi$ is a smooth manifold around $(0,0)$, with tangent space at such points $(w,z)$ given by
\[
T_{\mbox{\scriptsize gph\,} \Phi}\big(H(w),G(w,z)\big) ~=~
\big\{ \big(\nabla H(w)r , \nabla G(w,z)(r,s)\big) : r \in \W,~ s \in \Z \big\}.
\]
Its image under the projection map $\mbox{proj} \colon \gph \Phi \to \U$ is simply the range of $\nabla H(w)$.  Since
$\nabla H(0)$ is injective, the projection has locally constant rank $\dim \W$.  Partial smoothness follows, and the local description of the active manifold follows from Proposition \ref{active}.

Conversely, suppose $\Phi \colon \U \tto \V$ is partly smooth at $\bar u$ for $\bar v \in \Phi(\bar u)$.  By the Constant Rank Theorem, we can consider the projection map $\mbox{proj}$ as having the form $(w,z) \mapsto (w,0) \in \W \times \Y$, where $(w,z) \in \W \times \Z$ (for Euclidean spaces $\W$ and $\Z$) defines local coordinates for the manifold $\gph \Phi$, centered at $(\bar u,\bar v)$, and 
$(w,y) \in \W \times \Y$ (for a Euclidean space $\Y$) defines local coordinates for $\U$ centered around 
$\bar u$.

More explicitly, there exist maps
\begin{eqnarray*}
F \colon \W \times \Z \to \U, & \mbox{smooth around}~ (0,0), & \mbox{with}~ F(0,0) = \bar u \\
G \colon \W \times \Z \to \V, & \mbox{smooth around}~ (0,0), & \mbox{with}~ G(0,0) = \bar v \\
Q \colon \W \times \Y \to \U, & \mbox{smooth around}~ (0,0), & \mbox{with}~ Q(0,0) = \bar u
\end{eqnarray*}
with
\begin{eqnarray*}
\big(\nabla F(0,0), \nabla G(0,0)\big) \colon \W \times \Z & \to & \U \times \V \\
\nabla Q(0,0) \colon \W \times \Y & \to & \U
\end{eqnarray*}
both injective, and for all small $\delta > 0$,
\begin{eqnarray*}
\gph \Phi & = & 
\big\{ \big(F(w,z),G(w,z)\big) :  w \in B_{\delta}(0),~  z \in B_{\delta}(0) \big\} 
~~\mbox{around}~ (\bar u,\bar v)\\
\U & = & 
\big\{ Q(w,y) :  w \in B_{\delta}(0),~  y \in B_{\delta}(0) \big\}
~~\mbox{around}~ \bar v,
\end{eqnarray*}
and furthermore, $F(w,z) = Q(w,0)$ for all small $w \in \W$ and $z \in \Z$.

Now define a map $H \colon \W \to \U$ by $H(w) = Q(w,0)$, for $w \in \W$, and notice
$\nabla F(0,0) = (\nabla H(0),0)$.  Then, for points $w \in \W$ and $z \in \Z$, whenever
$0 = \nabla H(0)w = \nabla F(0,0)(w,z)$ and $\nabla G(0,0)(w,z) = 0$, we must have $w=0$ and $z=0$.  The result now follows.
\finpf

One consequence is the locally minimal identifiability of active manifolds we mentioned above, as we show next.

\begin{cor}[Minimal identifiability] \label{minimal}
If a set-valued mapping $\Phi \colon \U \tto \V$ is partly smooth at a point $\bar u \in \U$ for a value 
$\bar v \in \Phi(\bar u)$, then the active manifold $M$ has the following properties.
\begin{itemize}
\item
There exists a map $F \colon M \to \V$, smooth around $\bar u$, such that $F(\bar u) = \bar v$ and $F(u) \in \Phi(u)$ for all points $u \in M$ near $\bar u$.
\item
For any set $M' \subset \U$ containing $\bar u$, and any map $F' \colon M' \to \V$ that is continuous at 
'$\bar u$ and satisfies $F(\bar u) = \bar v$ and $F(u) \in \Phi(u)$ for all points $u \in M'$ near 
$\bar u$, we must have 
$M' \subset M$ around $\bar u$.
\item
$M$ is a locally minimal identifiable set at $\bar u$ for $\bar v$.
\end{itemize}
\end{cor}

\pf
To see the first property, we apply the coordinate representation guaranteed by Theorem \ref{coordinate}, and define the map $F$ by $F\big(H(w)\big) = G(w,0)$ for small vectors $w \in \W$.  The last property follows, since we just need to show the following inner semicontinuity property (see \cite[Proposition 2.8 ]{ident}:  for any sequence of points $u_r \to \bar u$ in the active manifold $M$, there exists a sequence of values $v_r \to \bar v$ with $v_r \in \Phi(u_r)$ for all large indices $r$.  To see this, simply set $v_r = F(u_r)$.

To see the second property, consider any sequence $u_r \in M'$ converging to $\bar u$.  By assumption, the sequence $\big(u_r,F'(u_r)\big) \in \gph \Phi$ converges to the point $(\bar u,\bar v)$, so $u_r \in M$ for all large indices $r$ by Proposition \ref{inclusion}.
\finpf

We also have the following calculus rule.

\begin{cor}[Sum rule] \label{Sum}
Consider a set-valued mapping $\Phi \colon \U \tto \V$ that is partly smooth at a point 
$\bar u \in \U$ for a value $\bar v \in \Phi(\bar u)$.  If the function $F \colon \U \to \V$ is smooth around $\bar u$, then the set-valued mapping $\Phi + F$ is partly smooth at $\bar u$ for 
$\bar v + F(\bar u)$, with the same dimension and active manifold.
\end{cor}

\pf
In terms of the coordinate representation guaranteed by Theorem \ref{coordinate}, we have
\[
\gph\!(\Phi + F) ~=~ 
\big\{ \big(H(w),\tilde G(w,z)\big) :  w \in B_{\delta}(0),~  z \in B_{\delta}(0) \big\}
~~\mbox{around}~ (\bar u,\bar v),
\]
where the map $\tilde G \colon \W \times \Z \to \V$ is defined by 
\[
\tilde G(w,z) = G(w,z) + F\big(H(w)\big)~~ (w \in \W,~ z \in \Z).  
\]
This map is smooth around the point $(0,0)$ with $\tilde G(0,0) = \bar v + F(\bar u)$.  Furthermore, by assumption,
\[
w \in \W,~ z \in \Z,~ \nabla H(0)w = 0~ \mbox{and}~ \nabla \tilde G(0,0)(w,z) = 0
~~\Rightarrow~~ w=0 ~\mbox{and}~ z=0,
\]
since
\[
\nabla \tilde G(0,0)(w,z) = \nabla G(0,0)(w,z) + \nabla F(\bar u)\nabla H(0)w.
\]
The result now follows by Theorem \ref{coordinate}.
\finpf

As with manifolds, a dual representation is sometimes more useful.

\begin{thm}[Dual representation] 
A set-valued mapping $\Phi \colon \U \tto \V$ is partly smooth at a point 
$\bar u \in \U$ for a value $\bar v \in \Phi(\bar u)$ if and only if it has a local representation of the following form:  there exist Euclidean spaces $\X$ and $\Y$, maps $P \colon \U \to \X$, smooth around 
$\bar u$ with $P(\bar u) = 0$ and $\nabla P(\bar u)$ surjective, and 
$Q \colon \U \times \V \to \Y$, smooth around $(\bar u,\bar v)$ with $Q(\bar u,\bar v) = 0$ and 
$\nabla_v Q(\bar u,\bar v)$ surjective, such that
\[
\gph \Phi ~=~ 
\big\{ (u,v) \in \U \times \V : P(u)=0,~  Q(u,v)=0 \big\}
~~\mbox{around}~ (\bar u,\bar v).
\]
The active manifold is then the inverse image $P^{-1}(0)$, around $\bar u$.
\end{thm}

\pf
Assuming the given representation, define a map $R \colon \U \times \V \to \X \times \Y$ by
$R(u,v) = \big( P(u),Q(u,v) \big)$ for points $u \in U$ and values $v \in V$.  Clearly $R$ is smooth around the point $(\bar u,\bar v)$, with $R(\bar u,\bar v) = (0,0)$.  The derivative
$\nabla R(\bar u,\bar v) \colon \U \times \V \to \X \times \Y$ is surjective, because for any values 
$x \in X$ and $y \in Y$ we can first find $r \in \U$ satisfying $\nabla P(\bar u)r = x$, and then find 
$s \in \V$ satisfying $\nabla_v Q(\bar u,\bar v)s = y - \nabla_u Q(\bar u,\bar v)r$, and in that case we have
\[
\nabla R(\bar u,\bar v)(r,s) ~=~ 
\big( \nabla P(\bar u)r, \nabla_u Q(\bar u,\bar v)r + \nabla_v Q(\bar u,\bar v)s \big)
~=~ (x,y).
\]
Since $\gph \Phi = R^{-1}(0,0)$ around the point $(\bar u,\bar v)$, we deduce that the graph of $\Phi$ is a manifold around $(\bar u,\bar v)$.  

For points $(u,v) \in \gph\Phi$ near the point $(\bar u,\bar v)$, we have
\begin{eqnarray*}
T_{\mbox{\scriptsize \gph} \Phi}(u,v) 
&=&
\mbox{Null}\big( \nabla R(u,v) \big) \\
&=&
\big\{ (r,s) \in \U \times  \V : \nabla P(u)r = 0,~\nabla_u Q(u,v)r + \nabla_v Q(u,v)s = 0 \big\},
\end{eqnarray*}
so, since the partial derivative $\nabla_v Q(u,v)$ is surjective, we deduce
\[
\mbox{proj} \big(T_{\mbox{\scriptsize \gph}\Phi}(u,v) \big) ~=~ \mbox{Null}\big( \nabla P(u) \big).
\]
Since the derivative $\nabla P(u)$ is surjective, this space has constant dimension for $u$ near $\bar u$, namely 
$\dim\U - \dim\X$, so partial smoothness follows.

Clearly the active manifold is contained in the inverse image $P^{-1}(0)$ around $\bar u$.  We claim these sets in fact agree around $\bar u$.  If not, there exists a sequence of points $u_k \to \bar u$ in $P^{-1}(0)$ lying outside the active manifold.  By the implicit function theorem, since the derivative 
$\nabla_v Q(\bar u,\bar v)$ is surjective, there exists a sequence of values $v_k \to \bar v$ such that $Q(u_k,v_k) = 0$ and hence $v_k \in \Phi(u_k)$ for all large $k$.  But this contradicts the definition of the active manifold.

Conversely,  suppose the mapping $\Phi$ is partly smooth at the point $\bar u \in \U$ for the value 
$\bar v \in \Phi(\bar u)$.  Using Theorem \ref{coordinate} (Coordinate representation), there exists a Euclidean space $\W$ and a map $H \colon \W \to \U$, smooth around $0$ with $H(0) = \bar u$ and derivative $\nabla H(0)$ injective, such that the active manifold is $M = H\big(B_{\delta}(0)\big)$ around 
$\bar u$ providing $\delta > 0$ is sufficiently small.  

Consider the map $G \colon \U \to \W$ discussed in Section \ref{Manifolds}, satisfying the property (\ref{inverse}), so its restriction $G|_M$ is the inverse of the diffeomorphism $H$ around the point $\bar u$.  Since 
$\gph \Phi$ is a manifold and contained in $M \times \V$ around the point $(\bar u,\bar v)$, the set
\[
\Lambda ~=~
\big\{ \big(G(u),v\big) : (u,v) \in \gph\Phi,~ u \in B_{\delta}(\bar u),~ v \in B_{\delta}(\bar v) \big\}
\]
is a manifold around the point $(0,\bar v) \in \W \times \V$.  Hence $\Lambda = S^{-1}(0)$ around 
$(0,\bar v)$, 
for some map $S \colon \W \times \V \to \Y$ (a Euclidean space), smooth around the point $(0,\bar v)$ with
$S(0,\bar v) = 0$ and $\nabla S(0,\bar v)$ surjective.  Equivalently, we have
\[
\gph \Phi ~=~ \big\{ \big(H(w),v\big) : S(w,v)=0,~w \in B_\delta(0),~v \in B_{\delta}(\bar v) \big\}~~
\mbox{around}~ (\bar u,\bar v).
\]

We claim, more precisely, that the partial derivative $\nabla_v S(0,\bar v) \colon \V \to \Y$ is surjective.  If not, there exists a nonzero vector $y \in \Y$ such that $\nabla_v S(0,\bar v)^* y = 0$.  By Corollary \ref{minimal} (Minimal identifiability), there exists a function $F \colon \W \to \V$, smooth around $0$, such that $F(0) = \bar v$ and $F(w) \in \Phi \big(H(w)\big)$ for all small vectors $w \in \W$.  We deduce $S\big(w,F(w)\big) = 0$ for all small $w \in \W$, so
\[
\nabla_w S(0,\bar v) + \nabla_v S(0,\bar v)\nabla F(0) = 0
\]
Taking adjoints shows $\nabla_w S(0,\bar v)^* y = 0$, so in fact $\nabla S(0,\bar v)^* y = 0$, contradicting the surjectivity of $\nabla S(0,\bar v)$.

There exists a Euclidean space $\X$ and a map $P \colon \U \to \X$, smooth around the point
$\bar u$, with $P(\bar u) = 0$ and $\nabla P(\bar u)$ surjective, such that the active manifold is 
$M = P^{-1}(0)$ around $\bar u$.  Furthermore, if we define a map $Q \colon \U \times \V \to \Y$ by
$Q(u,v) = S\big(G(u),v\big)$, then the desired representation now follows, since the partial derivative 
\[
\nabla_v Q(\bar u,\bar v) = \nabla_v S(0,\bar v)
\]
is surjective.
\finpf

\section{The normal bundle and partial smoothness}
Given a manifold $M \subset \U$ around a point $\bar u \in M$, we can consider the normal space as a set-valued mapping $N_M \colon \U \tto \U$, where we define $N_M(u) = \emptyset$ if $u \not\in M$.

\begin{thm}[Normal space] \label{normal}
If a set $M \subset \U$ is a $C^{(2)}$-smooth manifold around a point $\bar u \in M$, then the normal space mapping $N_M \colon \U \tto \U$ is partly smooth at $\bar u$ for any value $\bar v \in N_M(\bar u)$, with dimension $\dim\U$ and  active manifold $M$.
\end{thm}

\pf
We apply Theorem \ref{coordinate} (Coordinate representation).
Following the notation of Section \ref{Manifolds}, there exists a vector $\bar x \in \X$ satisfying 
$\nabla P (\bar u)^* \bar x = \bar v$.  We claim
\[
\gph N_M ~=~
\big\{ \big( H(w) , \nabla P \big(H(w)\big)^*x\big) : w \in B_{\delta}(0),~ x \in B_{\delta}(\bar x) \big\},~~
\mbox{around}~ (\bar u,\bar v),
\]
providing $\delta > 0$ is sufficiently small.  The inclusion ``$\supset$'' is clear, so it suffices to prove the inclusion ``$\subset$''.

For sufficiently small $\delta > 0$, the map $H$ gives a diffeomorphism between the open ball 
$B_\delta(0) \subset \W$ and an open neighborhood of the point $\bar u$ in the manifold $M$.  For such 
$\delta$, if the desired inclusion fails, then there exists a sequence of points $u_r \to \bar u$ in $M$ and a sequence of normals $v_r \in N_M(u_r)$ approaching $\bar v$, such that the sequence $(u_r,v_r)$ is disjoint from the right-hand side.  There must therefore exist a sequence of vectors $w_r \to 0$ in $\W$ satisfying $H(w_r) = u_r$, and a sequence of vectors $x_r \in \X$ satisfying 
\[
\nabla P(u_r)^*x_r = v_r \to \bar v = \nabla P (\bar u)^* \bar x.
\]
Since the linear map $\nabla P(\bar u)$ is surjective, we can represent it with respect to some orthonormal bases by the matrix $(A~ 0)$, where the matrix $A$ is invertible.  Denote the corresponding representation of $\nabla P(u_r)$ by $(A_r~ C_r)$, where $A_r \to A$ and $C_r \to 0$.  The property above ensures $A_r^T x_r \to A^T \bar x$ and hence $x_r \to \bar x$, contradicting our assumption that 
$x_r \not\in B_\delta(\bar x)$.

Now define a map
$G \colon \W \times \X \to \U$ by 
\[
G(w,z) ~=~ \nabla P \big(H(w)\big)^*(\bar x + z) ~~ \mbox{(for $w \in \W,~ z \in \X$)}.
\]
Clearly $G$ is smooth around the point $(0,0)$, with $G(0,0) = \bar v$.  Furthermore, around the point 
$(\bar u, \bar v)$, the graph of $\Phi$ has the representation (\ref{graph}), as we have just seen.  It remains to verify the regularity condition (\ref{regularity}).  By assumption, 
$\mbox{Null}\big(\nabla H(0)\big) = \{0\}$, so we just need to check that vectors $z \in \X$ satisfy the property
\[
\nabla G(0,0)(0,z) = 0 ~\Rightarrow~ z=0
\]
However, $\nabla G(0,0)(0,z) = \nabla P(\bar u)^* z$, and $\nabla P(\bar u)$ is surjective.  Notice that the dimension of $N_M$ is 
\[
\dim\W + \dim\X = \dim T_M(\bar u) + \dim N_M(\bar u) = \dim\U,
\]
so the result now follows.
\finpf

We can generalize this result substantially.  In the variational analysis that follows, we follow the terminology and notation of \cite{VA}.  The original definition of a partly smooth set appeared in \cite{Lewis-active}.  Here we use a slightly modified directional version \cite{ident}.

\begin{defn}
{\rm
Consider a closed set $Q \subset \U$, a point $\bar u \in Q$, and a normal vector 
$\bar v \in N_Q(\bar u)$.  We call $Q$ {\em partly smooth at $\bar u$ for $\bar v$ with respect to} a set $M \subset Q$ when all of the following properties hold.
\begin{itemize}
\item
$Q$ is prox-regular at $\bar u$ for $\bar v$.
\item
$M$ is a $C^{(2)}$-smooth manifold around $\bar u$.
\item
$N_M(\bar u) = \mbox{span}\,\hat N_Q(\bar u)$.
\item
For some neighborhood $W$ of $\bar v$, the mapping $u \mapsto N_Q(u) \cap W$ is inner semicontinuous at 
$\bar u$ relative to $M$.
\end{itemize}
}
\end{defn}

Since this definition is rather technical, a more concrete model is helpful.  Consider the {\em fully amenable} case when the set $Q$ coincides around $\bar u$ with an inverse image $F^{-1}(D)$ where $F$ is a 
$C^{(2)}$-smooth mapping and $D$ is a closed convex set satisfying $N_D\big(F(\bar u)\big) \cap N\big(\nabla F(\bar u)\big) = \{0\}$ (as holds in particular if $Q$ is closed and convex).  Then the prox-regularity condition holds, and the normal and regular normal cones, $N_Q(\bar u)$ and $\hat N_Q(\bar u)$, coincide.  The inner semicontinuity condition means that, for any normal vector $v \in N_Q(\bar u)$ near 
$\bar v$, and any sequence of points $u_r \to \bar u$ in $M$, there exists a corresponding sequence of normals $v_r \in N_Q(u_r)$ approaching $\bar v$.

We then have the following result.

\begin{thm}
Consider a closed set $Q \subset \U$, a point $\bar u \in Q$, a regular normal vector 
$\bar v \in \hat N_Q(\bar u)$, and suppose that $M \subset Q$ is a $C^{(2)}$-smooth manifold around 
$\bar u$.  Then the following properties are equivalent for the normal cone mapping $N_Q$.
\begin{enumerate}
\item[{\rm (i)}]
$N_Q$ is partly smooth at $\bar u$ for $\bar v$, with active manifold $M$.
\item[{\rm (ii)}]
$M$ is identifiable for $N_Q$ at $\bar u$ for $\bar v$. 
\item[{\rm (iii)}]
$Q$ is partly smooth at $\bar u$ for $\bar v$ with respect to $M$, and 
$\bar v \in \mbox{\rm ri}\,\hat N_Q(\bar u)$.
\item[{\rm (iv)}]
$\gph N_Q = \gph N_M$ around $(\bar u,\bar v)$.
\end{enumerate}
When these properties hold, the dimension of $N_M$ at $\bar u$ for $\bar v$ is just $\dim\U$.
\end{thm}

\pf
The implication (i) $\Rightarrow$ (ii) follows from Proposition \ref{inclusion}.  The equivalence of the properties (ii), (iii), and (iv) follows from \cite[Proposition 8.4]{ident}.  The implication (iv) $\Rightarrow$ (i) follows from Theorem \ref{normal}.
\finpf

The definition of a partly smooth function parallels that for sets.  Again we use a directional version of the original idea in \cite{Lewis-active}, following \cite{ident-arxiv}.

\begin{defn}
{\rm
Consider a closed function $f \colon \U \to \overline\R$, a point $\bar u \in \U$, and a subgradient 
$\bar v \in \partial f(\bar u)$.  We call $f$ {\em partly smooth at $\bar u$ for $\bar v$ with respect to} a set $M \subset \U$ when all of the following properties hold.
\begin{itemize}
\item
$f$ is prox-regular at $\bar u$ for $\bar v$.
\item
The restriction $f|_M$ is $C^{(2)}$-smooth around $\bar u$.
\item
The regular subdifferential $\hat \partial f(\bar u)$ is a translate of the normal space $N_M(\bar u)$.
\item
For some neighborhood $W$ of $\bar v$, the mapping $u \mapsto \partial f(u) \cap W$ is inner semicontinuous at $\bar u$ relative to $M$.
\end{itemize}
}
\end{defn}

We then have the following result.

\begin{thm} \label{psfunctions}
Consider a closed function $f \colon \U \to \overline\R$, a point $\bar u \in \U$, and a regular subgradient $\bar v \in \hat\partial f(\bar u)$.  Suppose that $f$ is subdifferentially continuous at 
$\bar u$ for $\bar v$.  Suppose furthermore that $M \subset Q$ is a $C^{(2)}$-smooth manifold around 
$\bar u$, and that the restriction $f|_M$ is $C^{(2)}$-smooth around $\bar u$.  Then there exists a function $\bar f \colon \U \to \R$ that is both $C^{(2)}$-smooth and satisfies $f|_M = \bar f|_M$ around 
$\bar u$, and for any such function  the following properties are equivalent for the subdifferential mapping $\partial f$.
\begin{enumerate}
\item[{\rm (i)}]
The mapping $\partial f$ is partly smooth at $\bar u$ for $\bar v$, with active manifold $M$.
\item[{\rm (ii)}]
The manifold $M$ is identifiable for $\partial f$ at $\bar u$ for $\bar v$. 
\item[{\rm (iii)}]
The function $f$ is partly smooth at $\bar u$ for $\bar v$ with respect to $M$, and 
$\bar v \in \mbox{\rm ri}\,\hat \partial f(\bar u)$.
\item[{\rm (iv)}]
Around $(\bar u,\bar v)$ we have
\[
\gph \partial f ~=~ \{(u,\nabla \bar f(u)+v) : u \in M,~v \in N_M(u)\}.
\]
\end{enumerate}
When these properties hold, the dimension of $\partial f$ at $\bar u$ for $\bar v$ is just $\dim\U$.
\end{thm}

\pf
The existence of the function $\bar f$ is just the definition smoothness of $f|_M$.
The implication (i) $\Rightarrow$ (ii) follows from Proposition \ref{inclusion}.  The equivalence of the properties (ii), (iii), and (iv) follows from \cite[Proposition 10.12]{ident-arxiv}.  The implication (iv) $\Rightarrow$ (i) follows from Theorem \ref{normal} and Corollary \ref{Sum} (Sum rule).
\finpf

Again, the assumptions are rather technical, so we illustrate with a more concrete model.  Consider the {\em fully amenable} case when the function $f$ is finite at $\bar u$ and agrees around $\bar u$ with a composite function $g \circ F$, where the mapping $F$ is $C^{(2)}$-smooth around $\bar u$ and the function $g$ is lower semicontinuous and convex, satisfying 
$N_{\mbox{\scriptsize \cl\!(\mbox{dom} \!g)}}\big(F(\bar u)\big) \cap N\big(\nabla F(\bar u)\big) 
= \{0\}$. (When $F$ is simply the identity mapping, we recover the case when $f$ is lower semicontinuous and convex).  Then both the subdifferential continuity and prox-regularity condition holds, and the normal and regular subdifferentials, $\partial f(\bar u)$ and $\hat \partial f(\bar u)$, coincide.

\section{Identifiability for primal-dual splitting}
We consider the saddlepoint problem
\[
\inf_{x \in \X} \sup_{y \in \Y} \{ (f+p)(x) + \ip{Ax}{y} - (g+q)(y) \}
\]
for Euclidean spaces $\X$ and $\Y$, lower-semicontinuous convex functions $f \colon \X \to \overline\R$ and $g \colon \Y \to \overline\R$, $C^{(2)}$-smooth convex functions $p \colon \X \to \R$ and $q \colon \Y \to \R$, and a linear map $A \colon \X \to \Y$.  Saddlepoints $(x,y)$ satisfy the inclusion
\[
(0,0) \in \Phi(x,y)
\]
where the set-valued mapping $\Phi \colon \X \times \Y \tto \X \times \Y$ is defined by
\[
\Phi(x,y) = 
\big(\partial f(x) + \nabla p(x) + L^* y\big) \times 
\big(-Lx + \partial g(y) + \nabla q(y)\big).
\]

The following method (following \cite{liang}) covers a variety of primal-dual algorithms \cite{chambolle-pock,vu,condat,he-yuan}.  As usual, we denote by $\mbox{prox}_f(x)$ the unique minimizer of the function 
$f(\cdot) + \frac{1}{2}\|\cdot-x\|^2$.
\newpage

\begin{alg}[Primal-dual splitting]
\label{shor-update-method}
{\rm
\begin{algorithmic}
\STATE
\STATE  Choose $\gamma,\mu > 0$.  For $k=0$, $x_0 \in \X$, $y_0 \in \Y$,
\WHILE{not done}
\STATE  $x_{k+1} = \mbox{prox}_{\gamma f}\big(x_k - \gamma \nabla p(x_k) - \gamma A^* y_k\big)$,
\STATE  $y_{k+1} = \mbox{prox}_{\mu g}\big(y_k - \mu \nabla q(y_k) + \mu A(2x_{k+1}-x_k)\big)$,
\STATE $k=k+1$;
\ENDWHILE
\end{algorithmic}
}
\end{alg}

\noindent
Assuming suitable conditions \cite[Theorem 3.3]{liang}, there exists a saddlepoint $(\bar x,\bar y)$ satisfying
\bmye \label{converges}
(x_k,y_k) \to (\bar x,\bar y) \qquad \mbox{and} \qquad
\mbox{dist}\big((0,0),\Phi(x_k,y_k)\big) \to 0.
\emye

Assume furthermore, again following \cite{liang}, that the function $f$ is partly smooth at $\bar x$ for 
$-\nabla p(\bar x) - L^* \bar y$ with respect to some set $M \subset \X$, that the function $g$ is partly smooth at $\bar y$ for $-\nabla q(\bar y) + L \bar x$ with respect to some set $N \subset \Y$, and that the nondegeneracy conditions
\[
-\nabla p(\bar x) - L^* \bar y \in \mbox{ri}\, \partial f(\bar x) \qquad \mbox{and} \qquad
-\nabla q(\bar y) + L \bar x \in \mbox{ri}\, \partial g(\bar y)
\]
hold.
Theorem \ref{psfunctions} implies that the mapping $\partial f$ is partly smooth at $\bar x$ for 
$-\nabla p(\bar x) - L^* \bar y$ with respect to $M$, and the mapping $\partial g$ is partly smooth at 
$\bar y$ for $-\nabla q(\bar y) + L \bar x$ with respect to $N$.  It follows immediately that the set-valued mapping $(x,y) \mapsto \partial f(x) \times \partial g(y)$ is partly smooth at $(\bar x,\bar y)$ for $(-\nabla p(\bar x) - L^* \bar y,-\nabla q(\bar y) + L \bar x)$ with respect to $M \times N$ and hence by the sum rule that the set-valued mapping $\Phi$  is partly smooth at $(\bar x,\bar y)$ for $(0,0)$ with respect to $M \times N$.  By Proposition \ref{inclusion}, $M \times N$ is identifiable for $\Phi$ at 
$(\bar x,\bar y)$ for $(0,0)$, so the convergence property (\ref{converges}) implies $x_k \in M$ and $y_k \in N$ eventually:  exactly the conclusion of \cite[Theorem 3.3]{liang}. 

\section{Example:  smooth optimization on a manifold}
We end with a brief but representative example to illustrate the interplay between partial smoothness and the second-order sufficient conditions.  
Suppose $M \subset \U$ is a $C^{(2)}$-smooth manifold around a point $\bar u \in M$, and 
$f \colon M \to \R$ is a $C^{(2)}$-smooth function.  We can 
consider a corresponding extended-valued function $\tilde f \colon \U \to \bR$ defined by
\[
\tilde f(u) = 
\left\{
\begin{array}{ll}
f(u)	& (u \in M) \\
+\infty	& (u \not\in M),
\end{array}
\right.
\] 
Its subdifferential map is given by
\[
\partial \tilde f(u) =
\left\{
\begin{array}{ll}
\nabla_M f(u) + N_M(u)	& (u \in M) \\
\emptyset				& (u \not\in M),
\end{array}
\right.
\]
where $\nabla_M f(u) \in T_M(u)$ denotes the covariant derivative.
By Corollary \ref{Sum} (Sum rule), this set-valued mapping $\partial \tilde f$ is partly smooth at 
$\bar u$ for any value in the set $\nabla_M f(\bar u) + N_M(\bar u)$.  In particular, assuming the first-order necessary condition 
\[
\nabla_M f(\bar u) = 0,
\]
then $\partial \tilde f$ is partly smooth at $\bar u$ for $0$, with dimension $\dim\U$ and active manifold $M$.

Now suppose further that $\bar u$ is a local minimizer around which $f$ grows quadratically:  for some $\delta > 0$,
\[
f(u) \ge f(\bar u) + \delta |u-\bar u|^2~~ \mbox{for all}~ u \in M ~\mbox{near}~ \bar u.
\]
Equivalently, in addition to the first-order condition, $f$ satisfies the second-order sufficient condition:  the covariant Hessian 
$\nabla_M^2 f(u) \colon T_M(u) \to T_M(u)$ (a self-adjoint linear map) is positive definite when 
$u = \bar u$.  We also have (from \cite{shanshan}):
\[
N_{\mbox{\scriptsize \gph}\,\partial \tilde f}(\bar u,0)
~=~
\big\{ (z,w) : w \in T_M(\bar u) ~\mbox{and}~ z + \nabla_M^2 f(\bar u)w \in N_M(\bar u) \big\}.
\]
Hence $\gph \partial \tilde f$ intersects the subspace 
$\U \times \{0\}$ transversally at $(\bar u,0)$.  To see this, note
\[
(z,w) \in N_{\mbox{\scriptsize \gph}\,\partial \tilde f}(\bar u,0) \cap N_{\U \times \{0\}}(\bar u,0)
\]
if and only if
\[
w \in T_M(\bar u),~~ z + \nabla_M^2 f(\bar u)w \in N_M(\bar u),~~ z=0.
\]
Since $\nabla_M^2 f(\bar u)$ is positive definite, the latter property holds if and only if $z=0$ and 
$w=0$, as required.  Consequently, $(\bar u,0)$ is an isolated transversal point of intersection of the two manifolds 
$\gph \partial \tilde f$ and $\U \times \{0\}$.

To summarize, satisfying the first-order optimality conditions for minimizing the smooth function $f$ on the manifold $M \subset \U$ amounts to finding a point in the intersection of the space $\U \times \{0\}$ and the graph of the subdifferential of the corresponding extended-valued function $\tilde f$.  Assuming the second-order sufficient conditions, the subdifferential is a partly smooth mapping of dimension 
$\dim \U$, and its graph intersects the subspace $\U \times \{0\}$ transversally at an isolated point.

\subsubsection*{Acknowledgements}
Many thanks to Artur Gorokh for many helpful comments during the development of these results.

\def\cprime{$'$} \def\cprime{$'$}


\end{document}